\documentclass[10pt]{amsart}
\usepackage {amssymb,latexsym}
\usepackage [all]{xy}

\theoremstyle{definition}

\numberwithin {equation}{section}

\def\cs{{$C^{\ast}$}}
\def\ra{{\rightarrow}}

\def\H{{\mathcal{H}}}
\def\gh{{(G,H)}}

\def\la{{\lambda}}

\def\inv{{^{-1}}}

\def\ba{{\backslash}}
\begin{document}
\title{An erratum to ``Property (RD) for Hecke pairs''}

\author{Vahid Shirbisheh}
\email{shirbisheh@gmail.com}




\maketitle


Let $\gh$ be a Hecke pair. The popular version of the reduced Hecke \cs-algebra associated to $\gh$ is based on the right action of $\H\gh$ (by the right convolution) on $\ell^2(G/H)$. In our paper \cite{s1}, since we wanted the Hecke algebra $\H\gh$ acts from the left on a Hilbert space, we used the the left action of $\H\gh$ (by the left convolution) on $\ell^2(H\ba G)$. To this work, we must have used the following formula for the convolution product instead of the one appearing in Definition 1.1 of \cite{s1}:
\begin{equation}
\label{1}
f_1\ast f_2 (g):=\sum_{\gamma \in <H\ba G>} f_1(g\gamma\inv)f_2 (\gamma),
\end{equation}
for all $f_1,f_2\in \H\gh$ and $g\in <G//H>$ (or similarly for all $f_1\in \H\gh$, $f_2\in \ell^2(H\ba G)$ and $g\in <H\ba G>$). One notes that these two formulas are equivalent for groups. But in the setting of Hecke pairs, our previous formula for the convolution product does not define a structure of an involutive algebra on $\H\gh$.

By applying this change to every occurrence of the convolution product in \cite{s1}, (precisely, in Pages 174, 182, 183), every statement is still valid.

To show that these two approaches for defining a reduced Hecke \cs-algebra are equivalent, let us denote the left and right actions of $\H\gh$ on $\ell^2(H\ba G)$ and $\ell^2(G/H)$ by $\la_l$ and $\la_r$, respectively.  Define $U:\ell^2(H\ba G)\ra \ell^2(G/H)$ by $U(\xi)(gH):=\xi(Hg\inv)$ for all $\xi\in \ell^2(H\ba G)$ and $g\in G$. It is clear that $U$ is a unitary operator from $\ell^2(H\ba G)$ onto $\ell^2(G/H)$. Using the fact that the map $\gamma\mapsto \gamma\inv$ is a bijection between the set $<H\ba G>$ of right cosets onto the set $<G/H>$ of left cosets, one easily checks $U(\la_l(f)(\xi))=\la_r(f)(U(\xi))$ for all $f\in \H\gh$ and $\xi\in \ell^2(H\ba G)$. In other words, these actions are unitarily equivalent, and so they define isomorphic concrete \cs-algebras acting on two different Hilbert spaces.

We would like to use this opportunity to correct a typo in the text of \cite{s1} as well. In Line 4 of Page 183 of \cite{s1} the summation should be over $G$ (instead of $<H\ba G>$).

\bibliographystyle {amsalpha}

\end{document}